\documentclass[12pt,reqno]{amsart} 

\usepackage[latin1]{inputenc}
\usepackage{amsmath, amsthm, amssymb,color,graphicx,bm,psfrag,float,verbatim}
\numberwithin{equation}{section}
\usepackage[all]{xy} 
\usepackage{fullpage,float} 
\newcommand{\ra}{\rightarrow} 
\newcommand{\la}{\leftarrow} 
\newcommand{\lra}{\longrightarrow} 
\renewcommand{\P}{\mathbb P} 
\newcommand{\conic}{\mathcal K}
\newcommand{\FF}{\mathfrak f} 
\newcommand{\RL}{\Theta} 
\newcommand{\RN}{\Phi} 
\newcommand{\CS}{\mathcal C} 
\newcommand{\fP}{\mathfrak P}
\newcommand{\F}{\mathcal F} 
\newcommand{\im}{\text{image}}
\newcommand{\complex}{\mathbf C} 
\newcommand{\six}{\text{\sc six}}
\newcommand{\SG}{\mathfrak S} 
\newcommand{\ltr}{\text{\sc ltr}}

\newcommand{\cZ}{\mathcal Z}
\newcommand{\bA}{\mathbb A} 
\newcommand{\bB}{\mathbb B} 
\newcommand{\bC}{\mathbb C} 
\newcommand{\bD}{\mathbb D} 
\newcommand{\bE}{\mathbb E} 
\newcommand{\bF}{\mathbb F}

\renewcommand{\proof}{{\sc Proof. \;}}

\newtheorem{Theorem}{Theorem}[section]
\newtheorem{Lemma}[Theorem]{Lemma}
\newtheorem{Proposition}[Theorem]{Proposition}

\newcommand{\pasc}[6]{\left\{ \begin{array}{ccc} #1 & #2 & #3\\ #4 &
     #5 & #6 \end{array} \right\}}

\begin{document} 
\title{When are the Cayley-Salmon lines conjugate?} 
\author{Jaydeep Chipalkatti} 
\maketitle

\bigskip 

\parbox{17cm}{ \small
{\sc Abstract:} Given six points on a conic, Pascal's
theorem gives rise to a well-known configuration called the
\emph{hexagrammum mysticum}. It consists of, amongst other things,
twenty Steiner points and twenty Cayley-Salmon lines. 
It is a classical theorem due to von Staudt that the Steiner points
fall into ten conjugate pairs with reference to the conic; but this is
not true of the C-S lines for a general choice of six points. 
It is shown in this paper that the C-S lines are pairwise conjugate precisely
when the original sextuple is~\emph{tri-involutive}. The variety of tri-involutive
sextuples turns out to be arithmetically Cohen-Macaulay of codimension
two. We determine its $SL_2$-equivariant minimal resolution.} 

\bigskip

AMS subject classification (2010): 14N05, 51N35. 

\bigskip 

\tableofcontents

\section{Introduction} \label{section.introduction} 

The main result of this paper involves the so-called Caylay-Salmon lines, which
form a subconfiguration of the famous~\emph{hexagrammum mysticum} of
Pascal. I have retained the notation of my earlier paper~\cite{CH1}, but some of the background
is reproduced below for ease of reading. 

\subsection{Pascal lines} 
Let $\conic$ denote a nonsingular conic in the complex projective
plane. Consider six distinct points $\Gamma = \{A, B, C, D, E, F \}$ on
$\conic$, arranged as an array $\left[ \begin{array}{ccc} A & B & C \\ F & E &
   D \end{array} \right]$. Pascal's theorem says that the three cross-hair intersection points
\[ AE \cap BF, \quad AD \cap CF, \quad BD \cap CE \] 
(corresponding to the three minors of the array) are collinear. The
line containing them (usually called the Pascal line, or just the
Pascal) is denoted by $\pasc{A}{B}{C}{F}{E}{D}$. It appears as
$k(1,23)$ in the usual labelling schema for Pascals (see~\cite[\S 2]{CH1}). 
Similarly, the lines $k(2,13)$ and $k(3,12)$ are respectively equal to 
\[ \pasc{A}{B}{C}{D}{F}{E} \quad \text{and} \quad 
\pasc{A}{B}{C}{E}{D}{F}. \] 

\thispagestyle{empty} 

\subsection{Steiner Points} 
It was proved by Steiner that the lines 
$k(1, 23), k(2, 13), k(3,12)$ are concurrent, and their common point of intersection is called a Steiner point,
denoted by\footnote{It is a notational convention that $123$ stands for the set
  $\{1,2,3\}$, and hence the order of indices is irrelevant. This 
  applies to all similar situations, so that $k(1,23)$ is the same as
  $k(1,32)$ etc.}
$G[123]$. This is understood to hold for any three indices
in the same pattern. If $\six$ denotes the set $\{1, 2, \dots, 6\}$, then for any three indices 
$x,y,z \in \six$, the Pascals 
\[ k(x,yz), \quad k(y,xz), \quad k(z,xy),  \] 
intersect in the point $G[xyz]$, giving altogether $\binom{6}{3} = 20$
Steiner points. This theorem, together with similar incidence theorems 
mentioned below may be found in Salmon's~\cite[Notes]{SalmonConics}, which also
contains references to earlier literature on the
subject. Salmon does not use the $k$-notation however, but there is a
detailed explanation of it in Baker's~\cite[Note II]{Baker}. 

\subsection{Kirkman points and Cayley-Salmon lines} 
It is a theorem due to Kirkman that the lines 
\[ k(1,23), \quad k(1,24), \quad k(1,34) \] 
are concurrent, and their common point of intersection is called the Kirkman
point $K[1,234]$. There are $60$ such points $K[w,xyz]$,
corresponding to indices $w \in \six, \{ x,y,z \} \subseteq \six
\setminus \{w\}$. It was proved by Cayley and
Salmon that the Kirkman points 
\[ K[4,123], \quad K[5,123], \quad K[6,123] \] 
are collinear, and the common line containing them is called the
Cayley-Salmon line\footnote{The illustrious 
$27$ lines on a nonsigular cubic surface
  (see~\cite[Ch.~4]{GH}) are also sometimes called the Cayley-Salmon
  lines. Although our context is different, there is, in fact, a 
  thematic connection between the geometry of the $27$ lines and the
  \emph{hexagrammum mysticum} (see~\cite{Richmond}).} $g(123)$. There are $20$ such lines $g(xyz)$
for $\{x,y,z\} \subseteq \six$. It is also the case that $g(xyz)$
contains the Steiner point $G[xyz]$. 

For a \emph{general} choice of six points $\Gamma$, all of the lines and
points above are pairwise distinct, and there are no further incidences
apart from those already mentioned. (It may happen that some of the lines and points become
undefined for special positions of $\Gamma$; this will be explained as
and when necessary.) Using the standard notation for
incidence correspondences (see~\cite{Grunbaum}), the situation for a
general $\Gamma$ can be summarised as follows: 
\begin{itemize} 
\item 
The $K$-points and $g$-lines make a $(60_1,20_3)$ configuration. That
is to say, each of the $60$ $K$-points lies on one $g$-line, and each
of the $20$ $g$-lines contains three $K$-points. 
\item 
The $K$-points and $k$-lines make a $(60_3,60_3)$ configuration. 
\item 
The $G$-points and $k$-lines make a $(20_3,60_1)$ configuration. 
\item 
The $G$-points and $g$-lines make a $(20_1,20_1)$ configuration. 
\end{itemize} 

We have throughout used uppercase letters for points, and lowercase
letters for lines. The governing pattern is that, if the $A$-points and $b$-lines form an $(r_s, t_u)$ 
configuration, then the $B$-points and $a$-lines form a $(t_u, r_s)$ configuration.\footnote{There are an
  additional fifteen $i$-lines and fifteen $I$-points which also obey this
  pattern, but we will not pursue this digression. See the notes by
  Baker and Salmon referred to above.} 
It is tempting to conjecture that this numerical duality should be
explainable as a pole-polar duality, but this is surprisingly not so (cf.~\cite[p.~194]{Hesse}). The following
pair of facts signals a breakdown in symmetry. 
Assume $\Gamma$ to be general, and let $\{u,v,w,x,y,z\} = \six$. Then 
\begin{itemize} 
\item 
The points $G[uvw]$ and $G[xyz]$ are conjugate with
respect to $\conic$; that is to say, the polar line of $G[uvw]$ passes
through $G[xyz]$ and conversely. This is a theorem due to von Staudt
(see~\cite{Staudt} or \cite[\S 86.2]{Pedoe}). 
\item 
The lines $g(uvw), g(xyz)$ are not conjugate; that is to say,
the pole of $g(uvw)$ does not lie on $g(xyz)$. 
This is easily checked with an example;
see~\S\ref{section.cs.nonconj} below. 
\end{itemize} 

It is natural to ask whether the $g$-lines are in fact pairwise
conjugate for \emph{special} positions of $\Gamma$. The main result of this
paper answers this question. 

\subsection{The tri-involutive configuration} \label{section.3inv.configuration}
Recall that $\conic$ is isomorphic to the projective line $\P^1$. 
We will say that a sextuple $\Gamma$ is \emph{tri-involutive}, if it
is projectively equivalent to 
\[ \left\{ 0, \, 1, \, \infty, \, p, \, \frac{p-1}{p}, \, \frac{1}{1-p} \right\}, \] 
for some $p \in \P^1$. (The rationale behind this name will be
clarified in \S\ref{section.3inv.geometry}. The configuration may well
have been classically known in some form, but it arose in
\cite[\S 4.5]{CH1} in the course of solving a different problem.) 

Since the points are assumed distinct, 
\[ p \neq 0, \, 1, \, \infty, \, \frac{1 \pm \sqrt{-3}}{2}. \] 
It may happen that the three Kirkman points $K[w,xyz], w \in \six \setminus \{x,y,z\}$
all coincide, so that the line $g(xyz)$ becomes undefined. 
A complete list of such cases is given in \S\ref{undefined.glines}. 

\subsection{} 
We will say that a sextuple $\Gamma$ satisfies Cayley-Salmon
conjugacy (CSC), if the following property holds: The lines
  $g(uvw)$ and $g(xyz)$ are conjugate, whenever both are defined and 
$\{u,v,w,x,y,z\} = \six$. 

\medskip 

This is the main result of the paper: 
\begin{Theorem} \rm 
The following statements are equivalent: 
\begin{enumerate} 
\item $\Gamma$ is tri-involutive. 
\item $\Gamma$ satisfies CSC. 
\end{enumerate} 
\label{main.theorem} \end{Theorem} 

Let $\Omega$ denote the Zariski closure of the subset of tri-involutive sextuples inside 
\[ \text{Sym}^6 \, \conic \simeq \text{Sym}^6 \, \P^1 \simeq \P^6. \] 
It is a four-dimensional irreducible projective subvariety. (The $p$ contributes one
parameter, and the additional three come from $SL_2$ acting on $\P^1$.) 
In \S\ref{section.3inv.group}, we will determine the group of symmetries of a generic tri-involutive sextuple
and use it to calculate the degree of $\Omega$. 
The main theorem is proved in \S\ref{section.cspoly}; it relies upon the properties of a collection of invariant
polynomials associated to $3$-element subsets of $\six$. These properties were discovered by explicitly calculating
and factoring these polynomials in {\sc Maple}. 

The $SL_2$-equivariant minimal resolution of the ideal of
$\Omega$ is determined in the concluding section. It implies that
$\Omega$ is an arithmetically Cohen-Macaulay subvariety of $\P^6$. 

\section{Preliminaries} \label{section.prelim}

\subsection{} The base field is throughout $\complex$. As in \cite[\S
3]{CH1}, let $S_d$ denote the vector space of homogeneous polynomials of
degree $d$ in the variables $\{x_1, x_2 \}$. We identify the projective plane with
$\P S_2$, and $\conic$ with the image of the Veronese imbedding $\P
S_1 \lra \P S_2$. A nonzero quadratic form in the $x_i$ simultaneously
represents a point in $\P^2$ and its polar line with respect to
$\conic$. All the line-intersections or point-joins can then be
calculated as transvectants of binary forms. In particular, two lines
represented by quadratic forms $\lambda$ and $\mu$ are conjugate 
exactly when the transvectant $(\lambda, \mu)_2$ is zero. 

The sextuple $\Gamma =
\{z_1, \dots, z_6\} \subseteq \complex \cup \{ \infty \} = \P^1$ is identified with the sextic form
$\prod\limits_{i=1}^6 (x_1 - z_i \, x_2)$. (If $z_i = \infty$, we
interpret the corresponding factor as $x_2$.) In this way, $\Gamma$
can be seen as a point of $\P S_6 \simeq \P^6$. 

\subsection{} \label{section.cs.nonconj} 
For instance, consider the sextuple 
\[ A = 0, \quad B = 1, \quad C = \infty, \quad D = -1, \quad E = 3,
\quad F = -5. \] 
Then $E$ is represented by $(x_1-3 \, x_2)^2$, and
$C$ by $x_2^2$ etc. Then the line $EC$ corresponds to 
\[ ((x_1-3 \, x_2)^2, x_2^2)_1 = x_2 \, (x_1 - 3 \, x_2), \] 
and one can calculate any required line or point from
\S\ref{section.introduction}. In particular, $g(123)$ is given by the
form $3 \, x_1^2 -26 \, x_1 \, x_2 -5 \, x_2^2$. All the other $g$-lines are calculated
similarly, and then it is straightforward to check that 
\[ (g(123), g(\alpha))_2 \neq 0, \] 
for any $3$-element subset $\alpha \subseteq \six$. That is to say, $g(123)$ is
not conjugate to any C-S line. Since CSC is a closed condition on
sextuples, it follows that a general sextuple does not satisfy CSC. 

\subsection{} \label{section.3inv.cr}  
Recall that the cross-ratio of four ordered points on $\P^1$ is defined to be 
\[ \langle z_1, z_2, z_3, z_4 \rangle =
\frac{(z_1-z_3) \, (z_2-z_4)}{(z_1-z_4) \, (z_2-z_3)}. \] 
If $r$ denotes the original cross-ratio, then permuting the
points in all possible ways gives six variants (see \cite[Ch.~1]{Seidenberg}), namely 
\[ r, \quad \frac{1}{r}, \quad 1-r, \quad \frac{1}{1-r}, \quad \frac{r-1}{r}, \quad
\frac{r}{r-1}. \] 

One can reformulate the notion of tri-involutivity in terms of
cross-ratios. As in \cite[\S 3]{CH1}, consider the set of
letters $\ltr = \{\bA, \bB, \bC, \bD, \bE, \bF \}$, and define a \emph{hexad} to be 
an injective map $\ltr \stackrel{h}{\lra} \conic$. 
Write $h(\bA) =A, h(\bB) = B$, and so on for the corresponding points on
$\conic$. Then a sextuple $\Gamma$ is 
tri-involutive, if and only if there exists a hexad with image $\Gamma$ such that 
\begin{equation} \langle A,B,C,F \rangle = \langle B, C,
A, E \rangle =  \langle C,A,B,D \rangle. 
\label{crossratios.eq1} \end{equation} 
Indeed, if we let $A = 0, B = 1, C = \infty$, and if $p$ denotes this common
cross-ratio, then 
\begin{equation} D = p, \quad E = \frac{p-1}{p}, \quad F =
  \frac{1}{1-p}. 
\label{3inv.stdvalues} \end{equation} 

\subsection{} \label{section.3inv.geometry}
\begin{figure}
\includegraphics[width=8cm]{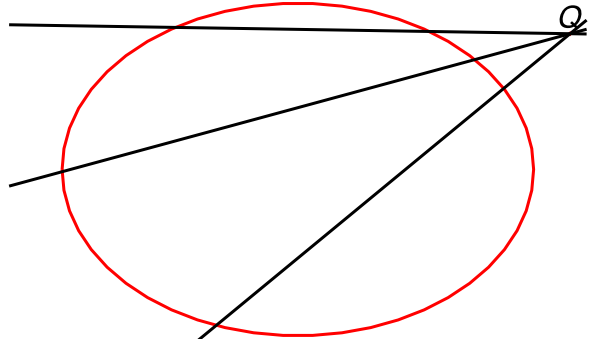}
\caption{A sextuple in involution} 
\label{fig:inv} 
\end{figure}

Recall (\cite[\S3]{CH1}) that a sextuple $\Gamma$ is said to be in
involution (see Figure~\ref{fig:inv}), if there exists a point $Q \in \P^2 \setminus \conic$ and
three lines $L, L', L''$ through $Q$ such that $\Gamma = \conic \cap (L \cup L' \cup L'')$. 
One says that $Q$ is a centre of involution for $\Gamma$. 

\begin{figure}
\includegraphics[width=12cm]{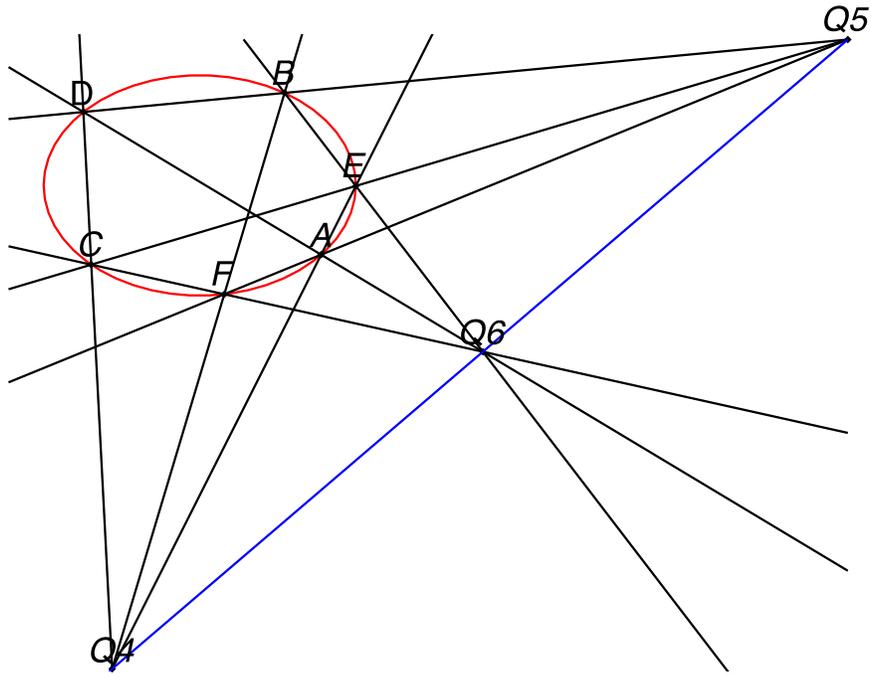}
\caption{A tri-involutive sextuple}
\label{fig:3inv} 
\end{figure}

Now a tri-involutive sextuple has the peculiarity that it has three
such centres, which are in fact collinear. (This is
explained at length in \cite[\S 4.5]{CH1}.) Observe that in
Figure~\ref{fig:3inv} on the next page, each of the three points $Q_4, Q_5, Q_6$ is a
centre of involution in such a way that 
\begin{itemize} 
\item $AE, CD, BF$ intersect in $Q_4$,  
\item $AF, CE, BD$ intersect in $Q_5$,   
\item $AD, BE, CF$ intersect in $Q_6$. 
\end{itemize}

\subsection{The exotic isomorphism} \label{section.exoticiso} 
It is well-known that the permutation group on six objects has a unique
outer automorphism (see \cite{CoxeterS6} or~\cite{Howard_etal}). Part of its charm is that it is implicated
everywhere in the geometry of the \emph{hexagrammum mysticum}.  We
will use it in the following version: let $\SG(X)$ denote the
permutation group on the set $X$. 
Then the following symmetric table defines an isomorphism 
$\SG(\ltr) \stackrel{\zeta}{\lra} \SG(\six)$. 

\[ 
\begin{array}{|c|c|c|c|c|c|c|} \hline 
{} & \bA & \bB & \bC & \bD & \bE & \bF \\ \hline 
\bA & {} & 14.25.36 & 16.24.35 & 13.26.45 & 12.34.56 & 15.23.46 \\ 
\bB & 14.25.36 & {} & 15.26.34 & 12.35.46 & 16.23.45 & 13.24.56 \\ 
\bC & 16.24.35 & 15.26.34 & {} & 14.23.56 & 13.25.46 & 12.36.45 \\ 
\bD & 13.26.45 & 12.35.46 & 14.23.56 & {} & 15.24.36 & 16.25.34 \\ 
\bE & 12.34.56 & 16.23.45 & 13.25.46 & 15.24.36 & {} & 14.26.35 \\ 
\bF & 15.23.46 & 13.24.56 & 12.36.45 & 16.25.34 & 14.26.35 & {} \\ \hline 
\end{array} \]

For instance, the entry in row $\bA$ and column $\bB$ is $14.25.36$,
and hence $\zeta$ takes the transposition $(\bA \, \bB)$ to the element 
$(1 \, 4) \, (2 \, 5) \, (3 \, 6)$ of cycle type
$2+2+2$. In the reverse direction, $\zeta^{-1}$ takes $(1 \, 2)$ to 
$(\bA \, \bE) \, (\bB \, \bD) \, (\bC \, \bF)$, because these are
precisely the positions in the table where the pair $12$ appears. The
action of $\zeta$ or $\zeta^{-1}$ on an arbitrary element can be read off by
writing it as a product of transpositions. 

\section{The symmetries of a tri-involutive sextuple} 
\label{section.3inv.group} 

\subsection{}  Let $\Gamma$ be a tri-involutive sextuple. We will say that a hexad
$\ltr \stackrel{h}{\lra} \Gamma$ is an \emph{alignment} if equations
(\ref{crossratios.eq1}) hold. It would be convenient to characterise
all shuffles of $\ltr$ which preserve this property. After an
automorphism of $\P^1$, we may assume that $h$ is given by 
\[ \bA \ra 0, \quad \bB \ra 1, \quad \bC \ra \infty, \quad 
\bD \ra p, \quad \bE \ra \frac{p-1}{p}, \quad \bF \ra
\frac{1}{1-p}. \] 
Let $\RL(p)$ denote the subgroup of bijections $\ltr \stackrel{\eta}{\lra} \ltr$ such
that 
\begin{equation} \langle A', B', C', F' \rangle = 
\langle B', C', A', E' \rangle =  
\langle C', A', B', D' \rangle, 
\label{crossratios.eq3} \end{equation} 
where $A'$ stands for $h \circ \eta(\bA)$ and
similarly for $B'$ etc. The group is completely determined by formal
properties of cross-ratios for a general value of $p$, and as such is independent of $p$. We 
will denote the generic case simply by $\RL$, and proceed to determine
its structure. (The group may be larger for special values of $p$;
see~\S\ref{gorenstein3.case} below.) This will be used later to find the degree of the locus of
tri-involutive sextuples. 

With the alignment as above, 
\begin{equation} \langle E,F,D,B \rangle = \langle F, D, 
E, A \rangle =  \langle D,E,F,C \rangle = p, 
\label{crossratios.eq2} \end{equation}  
and hence $(\bA \, \bE) \, (\bB \, \bF) \, (\bC \, \bD) \in \RL$. More generally, it is easy to verify that 
\begin{equation} (x \, x') \, (y \, y') \, (z \,
 z') \in \RL, \label{RL.elements} \end{equation} 
whenever 
\[ \{x,y,z\} = \{\bA,\bB,\bC\}, \quad \text{and} \quad \{x',y',z'\} =
\{\bD,\bE,\bF\}. \] 
(For instance, if $\eta = (\bA \, \bD) \, (\bB \, \bF) \, (\bC \, \bE)$, then each
term in (\ref{crossratios.eq3}) becomes $p/(p-1)$.) It immediately
follows that $(x \, x') \, (y \, y') \in \RL$ whenever $x,x' \in
\{\bA,\bB,\bC\}$ and $y,y' \in \{\bD,\bE,\bF\}$. Moreover, since 
$(\bA \, \bB) \, (\bD \, \bE) \, (\bA \, \bC) \, (\bD \, \bE) = (\bA
\, \bB \, \bC)$ etc, every $3$-cycle in either $\SG(\{\bA,\bB,\bC\})$ or $\SG(\{\bD,\bE,\bF\})$
is in $\RL$.

\begin{Proposition} \rm 
For a general value of $p$, the elements in (\ref{RL.elements}) generate $\RL$. 
\end{Proposition} 
\proof 
Let $\RL' \subseteq \RL$ denote the subgroup they generate, and
consider the subgroup 
\begin{equation} \RN = \SG(\{1,2,3\}) \times \SG(\{4,5,6\}), 
\label{defn.RN} \end{equation}
of $\SG(\six)$. We claim that there is an equality 
$\zeta(\RL') = \RN$. Indeed, the table in \S\ref{section.exoticiso}
shows directly that $\zeta$ takes the six
elements in (\ref{RL.elements}) to the $2$-cycles 
$(x \, x')$, where $x,x'$ are both in $\{1,2,3\}$ or in
$\{4,5,6\}$. (For instance, $(\bA \, \bE) \, (\bB \, \bF) \, (\bC \,
\bD)$ goes to $(5 \, 6)$.) This proves the claim. 

Now assume that $\RL' \subsetneq \RL$, and let $z \in \zeta(\RL) \setminus \RN$ denote an element
with the maximum number of fixed points. It must be the case that $z$
takes elements in $\{1,2,3\}$ to $\{4,5,6\}$ and conversely; for if
not, one could increase the number of fixed points by multiplying it
by an element of $\RN$.  But then after conjugation by an element of
$\RN$, one may assume $z$ to be one of the following: 
$(1 \, 4), (1 \, 4) \, (2 \, 5), (1 \, 4) \, (2 \, 5) \, (3 \, 6), 
(1 \, 4 \, 2 \, 5), (1 \, 4 \, 2 \, 5) \, (3 \, 6), (1 \, 4 \, 2 \, 5
\, 3 \, 6)$. Now apply $\zeta^{-1}$ to each of these, and
check that none of the images belongs to $\RL$. This proves
that $\RL' = \RL$. \qed 

It follows that $|\RL| =36$.  The isomorphism $\RL \simeq \RN$ shows a common phenomenon
surrounding the exotic isomorphism; 
namely a complicated structure on one side of $\zeta$ often appears as a
simpler structure on the other side. 

\subsection{} 
Now let $\Omega \subseteq \P^6$ denote the Zariski closure of the set of 
tri-involutive sextuples. It is an irreducible projective fourfold. 

\begin{Proposition} \rm 
The degree of $\Omega$ is $16$. 
\label{prop.degree16} \end{Proposition} 

The proof will emerge from the discussion below. 
If $z \in \conic$ is an arbitrary point, then $\{ \Gamma \in
\P^6: z \in \Gamma \}$ is a hyperplane in $\P^6$. Since the
degree of $\Omega$ is the number of points in its 
intersection with four general hyperplanes, we are reduced to the
following question: Given a set of four general points $Z = \{z_1, \dots, z_4
\} \subseteq \conic \simeq \P^1$, find the number of tri-involutive sextuples $\Gamma$ such that 
$Z \subseteq \Gamma$. In other words, find all tri-involutive extensions of $Z$. 

Let the variables $a,b, \dots, f$ respectively stand for the
coordinates of $A, B, \dots, F$ on $\conic$. Then conditions~(\ref{crossratios.eq1})
correspond to equations 
\begin{equation} \frac{(a-c)(b-f)}{(a-f)(b-c)} = \frac{(b-a)(c-e)}{(b-e)(c-a)} =
\frac{(c-b)(a-d)}{(c-d)(a-b)}. 
\label{crdeg.eq} \end{equation}

\subsection{} 
By an \emph{assignment}, we will mean an injective map $Z \stackrel{u}{\lra} \{a, \dots,
f\}$. Given such an assignment, say 
\[ z_1 \ra c, \quad z_2 \ra e, \quad z_3 \ra a, \quad z_4 \ra d, \] 
we abbreviate it by $[c,e,a,d]$. Now substitute $c = z_1, \dots, d
= z_4$ into (\ref{crdeg.eq}), and solve 
for the remaining unknowns $b$ and $f$. Depending on the number of
solutions so obtained, one would get one or more tri-involutive
sextuples $\Gamma$ containing $Z$. It it clear that every extension must come from an assignment. 
The group $\RL$ acts naturally on the set of assignments by
permuting the $a, \dots, f$, and two assignments which are in the
same orbit have the same extensions. 
Given $u$ as above, consider the integer $n(u) = \text{card} \, (\im(u) \cap \{a,b,c\})$. 
It must be either $1,2$ or $3$. 

\subsection{} Assume $n(u)=3$. After an action by an element of $\RL$,
we may assume that $u$ is one of the following: 
\begin{equation} [d,a,b,c], \quad [a,d,b,c], \quad [a,b,d,c], \quad
  [a,b,c,d]. 
\label{ext.4choices} \end{equation} 
(For example, $[c,a,f,b]$ can be changed to $[a,c,e,b]$ by 
$(\bA \, \bC) \, (\bE \, \bF)$, and then to $[a,b,d,c]$ by $(\bB \, \bC) \, (\bD \, \bE)$.) In each case,
(\ref{crdeg.eq}) give linear equations for the 
remaining variables $e$ and $f$ which determine them uniquely. For instance, if we let 
the $z_i$ to be respectively $0, 1,\infty, r$, then the $[d,a,b,c]$ case gives the extension 
$Z \cup \left\{ \frac{r^2-r+1}{r}, r^2-r+1\right\}$. The remaining cases in (\ref{ext.4choices}) give
three more extensions by the following enlargements of $Z$: 
\[ \left\{ \frac{r^2}{r-1}, r-r^2 \right\}, \quad 
\left\{\frac{r}{r^2-r+1}, \frac{r^2}{r^2-r+1} \right\}, \quad 
\left\{ \frac{1}{1-r}, \frac{r-1}{r} \right\}. \] 
Thus altogether we get four extensions. The assumption $n(u)=1$ gives nothing
new, since we are back to $n(u)=3$ after an action by $(\bA \, \bD) \, (\bB \, \bE) \, (\bC \,
\bF)$. 

\subsection{} If $n(u)=2$, then after an action of $\RL$ we may assume that $\im(u) =
\{a,b,e,f\}$. We can then reduce further to one of the following six cases: 
\[ [a,b,x,y], \quad [a,x,b,y], \quad [a,x,y,b], \quad \text{where} \; 
\{x,y\} = \{e,f\}. \] 
(For instance, $[f,b,a,e]$ is changed to $[b,f,e,a]$ by $(\bA \, \bE)
\, (\bB \, \bF) (\bC \, \bD)$, and then to $[a,e,f,b]$ by $(\bA \,
\bB) \, (\bE \, \bF)$.) Given any such $u$, the first pair of
equations in (\ref{crdeg.eq}) gives a quadratic in $c$, and either of its
roots will determine the value of $d$ uniquely. Thus we get $6 \times
2 = 12$ extensions altogether, and the fact
that they are generically distinct is verified by a direct calculation
as above. Thus $\deg \Omega = 4+ 12 = 16$. \qed

\section{The Cayley-Salmon polynomials} \label{section.cspoly} 

Consider the graded polynomial ring $\fP =
\complex[a,b,c,d,e,f]$. One can formulate the CSC condition as the vanishing of certain
homogeneous polynomials in $\fP$; this will lead to a proof of the main theorem. 

\subsection{} 
As indicated earlier, each of the lines and points mentioned in
\S\ref{section.introduction} is represented by a quadratic form in
$x_1$ and $x_2$, with coefficients in $\fP$. For
instance, the line $k(1,23)$ appears as $\frac{1}{4}(c-d)$  times 
\[ (c \, f-c \, e-a \, d-b \, f+b \, d+a \, e) \, x_1^2 + \dots + 
(a \, c \, e \, f-a \, c \, d \, e+b \, c \, d \, f-b \, c \, e \, f+a
\, b \, d \, e-a \, b \, d \, f) \, x_2^2. \] 
Since $c - d \neq 0$, we are free to ignore the extraneous multiplicative
factor, and this will always be done in the sequel. Henceforth we will
not write out such forms explicitly, since the expressions are
generally lengthy and not much is to be learned by merely looking at them. I have programmed
the entire procedure in {\sc Maple}, and the structure of the polynomials
$\CS_\alpha$ described below was discovered in this way. 

\subsection{} If quadratic forms $Q, Q'$ represent the lines $g(123), g(456)$, then the vanishing of 
$(Q,Q')_2$ is the necessary and sufficient condition for the
lines to be conjugate. We denote this expression (shorn of
extraneous multiplicative factors) by $\CS_{123}$ or $\CS_{456}$; and
similarly define the \emph{Cayley-Salmon polynomial} $\CS_\alpha \in \fP$ for
each $3$-element set $\alpha \subseteq \six$. (Thus, $C_\alpha = C_{\six \setminus \alpha}$ by construction.) It turns out that
each $\CS_\alpha$ is homogeneous of degree $18$; moreover,
it is a product of six irreducible factors which
exhibit a high degree of symmetry. We will describe the structure of
$\CS_{123}$ in detail, and then the action of the
permutation group can be used to infer it for any $\CS_\alpha$. 

An explicit factorisation (done in {\sc Maple}) shows that 
\[ \CS_{123} = L_1  \, L_2 \, L_3 \, M_4 \, M_5 \, M_6, \] 
where each $L_i$ and $M_j$ is a homogeneous cubic in
$a, \dots, f$. (The factors are determined only up to nonzero 
multiplicative scalars, but this will cause no difficulty.) For instance, 
\begin{equation} 
L_1 = (c \, d \, e +b \, c \, d+a \, d \, f+a \, b \, e+a \, b \, f+b
\, e \, f) - (a \, c \, d+b \, d \, e+c \, d \, f+b \, c \, e+a \, c
\, f+a \, e \, f).  
\label{L1.expression} \end{equation} 
The following discussion will clarify how the remaining factors are thereby
determined, and why they are so labelled. 

\subsection{} 
The group $\SG(\ltr)$ acts naturally on $\fP$, and so does the group $\SG(\six)$ via
$\zeta^{-1}$. Let $G = G_{123} = G_{456}$ denote the subgroup of
$\SG(\six)$ generated by $\SG(\{1,2,3\})$ and $\SG(\{4,5,6\})$,
together with all elements of the form $(x \, x') \, (y \, y') \, (z \, z')$, 
where $\{x,y,z\} = \{1,2,3\}$ and $\{x',y',z'\} = \{4,5,6\}$. In other
words, $G$ is the subset of elements in $\SG(\six)$ which stabilise the
unordered pair of sets $\{\{1,2,3\}, \{4,5,6\}\}$. 

A direct calculation shows that $G$ permutes the set of factors $\F_{123} = \{L_1, L_2, L_3, M_4, M_5, M_6
\}$; specifically, we have a homomorphism $G \lra \SG(\F_{123})$ given by 
\[ \begin{array}{ll} 
(1 \, 2) \lra (L_1 \, L_2), & (1 \, 3) \lra (L_1 \, L_3), \\ 
(4 \, 5) \lra (M_4 \, M_5), & (4 \, 6) \lra (M_4 \, M_6), \\ 
(1 \, 4) \, (2 \, 5) \, (3 \, 6) \lra (L_1 \, M_4) \, (L_2 \, M_5) \,
(L_3 \, M_6). \end{array} \] 
This is interpreted as follows. Since $\zeta^{-1}(1 \, 2) = (\bA \, \bE) \,
(\bB \, \bD) \, (\bC \, \bF)$, the simultaneous substitution 
\[ a \ra e, \quad e \ra a, \quad b \ra d, 
\quad d \ra b, \quad c \ra f, \quad f \ra c, \] 
interchanges $L_1$ and $L_2$. The same recipe
applies throughout; for instance, since 
\[ \zeta^{-1}((1 \, 4) \, (2 \, 5) \, (3 \, 6)) = (\bA \, \bB), \] 
the simultaneous substitution $a \ra b, b \ra a$ acts as 
the permutation $(L_1 \, M_4) \, (L_2 \, M_5) \, (L_3 \, M_6)$. 

The $L_i$ and $M_j$ are labelled so as to harmonise with the action 
of $G$.  Since this action is transitive on $\F_{123}$, each of the
remaining factors $L_2, \dots, M_6$ is obtainable by a change of
variables in $L_1$. 

\subsection{Invariance} Consider a fractional linear transformation
(FLT) \[ \mu(z) = \frac{p \, z + q}{r\, z + s}, \] 
where $p, q, r, s$ and $z$ are, at the moment, seen as formal indeterminates. 
Then a calculation shows that 
\[ L_1(\mu(a), \dots, \mu(f))  = \frac{(p \, s - q \, r)^3} {(r \, a + s)
  \dots (r \, f + s)} \; L_1(a,\dots, f). \] 
In other words, $L_1$ is an invariant for the action of FLTs on ordered sextuples in
$\P^1$. A general theory of such functions 
may be found in~\cite[\S 1]{DolgachevOrtland}. 

\subsection{} 
Given a hexad $\ltr \stackrel{h}{\lra} \Gamma \subseteq \complex$, one can evaluate each
of the polynomials above by substituting $a = h(\bA), \dots, f = h(\bF)$. We will disallow
the value $\infty$ for simplicity, although it would not have been difficult
to account for this possibility. 

\begin{Proposition} \rm 
Assume that $h$ is an alignment. Then $M_4 = M_5 = M_6 =0$. 
\end{Proposition} 
\proof 
By the invariance mentioned above, the vanishing is unaffected by an
FLT of $\Gamma$, hence we may assume $a=1,b=0,c=-1$. If $t$ denotes
the common cross-ratio in (\ref{crossratios.eq1}), then we must have 
\[ d = \frac{1-t}{1+t}, \quad e = \frac{1}{2 \, t -1}, \quad f = \frac{t}{t-2}. \] 
Make a change of variables $a \ra b, b \ra a$ into $L_1$ and substitute the
values above; then the expression is seen to vanish. This proves that
$M_4=0$. Now recall that $\zeta^{-1}(4 \, 5) \in \RL$, i.e., 
the substitution $\zeta^{-1}(4 \, 5) = (\bA \, \bD) \, (\bB \, \bE) \,
(\bC \, \bF)$ followed by $h$ remains an alignment. This shows that $M_5=0$, and likewise for $M_6$. \qed 

\subsection{} 
Now assume $\alpha = \{1,2,4\}$. We have a factorisation 
\[ \CS_{124} = L_1'  \, L_2' \, L_4' \, M_3' \, M_5' \, M_6', \] 
with similar properties. Define the group $G_{124}$ analogously, and let $\F_{124} = \{L_1', L_2', L_4', M_3', M_5',
M_6'\}$. As before, we have a homomorphism $G_{124} \lra \SG(\F_{124})$ given by 
\[ \begin{array}{ll} 
(1 \, 2) \lra (L_1' \, L_2'), & (1 \, 4) \lra (L_1' \, L_4'), \\ 
(3 \, 5) \lra (M_3' \, M_5'), & (3 \, 6) \lra (M_3' \, M_6'), \\ 
(1 \, 3) \, (2 \, 5) \, (4 \, 6) \lra (L_1' \, M_3') \, (L_2' \, M_5') \,
(L_4' \, M_6'). \end{array} \] 
The element $(3 \, 4)$ takes the set $\{1,2,3\}$ to $\{1,2,4\}$, and 
the change of variables induced by 
\[ (3 \, 4) \lra (\bA \, \bE) \, (\bB \, \bC) \, (\bD \, \bF) \] 
gives a bijection $\beta_{(3,4)}: \F_{123} \lra \F_{124}$ as follows: 
\[ \begin{array}{lll} 
L_1 \lra L_1', & L_2 \lra L_2', & L_3 \lra L_4', \\ 
M_4 \lra M_3', & M_5 \lra M_5', & M_6 \lra M_6'. 
\end{array} \] 
In this way, the structure of each $\F_\alpha$ is determined by that
of $\F_{123}$. It follows that the invariance property is also true of
each of the factors of $\CS_\alpha$. 
\begin{Lemma} \rm 
The forms $M_4$ and $M_3'$ are equal (up to a scalar). 
\end{Lemma} 
\proof 
Let $u = (1 \, 4) \, (2 \, 5) \, (3 \, 6)$ and $v = (3 \, 4)$. Since $u
(L_1) = M_4$ and $v(M_4) = M_3'$, it would suffice to show that
$u^{-1} \, v \, u = (1 \, 6)$ stabilises $L_1$ up to a scalar. But $(1 \, 6) \lra 
(\bA \, \bC) \, (\bB \, \bE) \, (\bD \, \bF)$, and it is easy to see
that the corresponding change of variables turns $L_1$ into $-L_1$. \qed

In particular, if $h$ is an alignment, then $M_3'$ vanishes and hence so
does $\CS_{124}$. 

\begin{Proposition} \rm 
Assume that $\Gamma$ is a tri-involutive sextuple. Then $\Gamma$
satisfies CSC. 
\end{Proposition} 
\proof 
Fix an alignment $h$, and consider a $3$-element subset $\alpha \subseteq \six$. 
By changing $\alpha$ to $\six \setminus \alpha$ if necessary, we may
assume that $\alpha$ has at least two elements in common with $\{1,2,3\}$. But then we are in the situation
above up to a change of indices, hence $\CS_\alpha=0$. (This is best
seen with an example. Say, $\alpha = \{2,4,6\}$, which we may change
to $\{1,3,5\}$. Now the element $\zeta^{-1} ((2 \, 3) \, (4 \, 5)) \in
\Theta$ preserves $\F_{123}$ and induces a bijection between
$\F_{124}$ and $\F_{135}$. Hence $\CS_{135}=0$.) \qed 

This proves the implication $(1) \Rightarrow (2)$ of the main theorem. 

\subsection{} 
For the converse, assume that $\ltr \stackrel{h}{\lra} \Gamma$ is a
hexad such that each $\CS_\alpha$ is zero. Then (at least) one of the factors in $\F_{123}$
must vanish, and we may assume it to be $L_1$ after a change of
variables. Let $a=1,b=0,c=-1$ after an FLT, then~(\ref{L1.expression}) 
reduces to 
\begin{equation} 
 L_1 = d \, (1 - e + 2 \, f) - f \, (1+e) = 0.
\label{L1.eq} \end{equation} 

Now (at least) one of the six factors in $\F_{124}$ must also vanish. 
When we pair each of them in turn with $L_1$, there are six cases
to consider. First, assume that 
\begin{equation} L_1' = d + e - f - d \, e \, f =0. 
\label{L1prime.eq} \end{equation} 
Now these two equations have a general solution 
\[ d = t, \quad e = \frac{2 \, t^2}{1 + t^2}, \quad f  = \frac{t \, (t+1)}{2 \, t^2
  - t + 1};   \] 
for a parameter $t$. But then 
\[ \bA \ra 1, \quad \bB \ra 0, \quad \bC \ra t, \quad 
\bD \ra -1, \quad \bE \ra \frac{t \, (t+1)}{2 \, t^2- t + 1}, \quad \bF \ra
\frac{2 \, t^2}{1+t^2} \] 
is an alignment, and hence $\Gamma$ must be tri-involutive. 
Replacing $L_1'$ with $L_4', M_5'$ or $M_6'$ leads to similar
solutions, and to the same conclusion. 

\subsection{} The remaining two cases behave a little differently. 
Assume that $L_1=0$ as above, and 
\[ L_2' = d-d \, e-d \, f+2 \, e \, f-d \, e \, f =0. \] 
The general solution is 
\[ d = \frac{t}{2+t}, \quad e = \frac{t-1}{t+3}, \quad f = t, \] 
which by itself does not force $\Gamma$ to be tri-involutive. However,
substituting it into $\CS_{125}$ leads to the equation 
\[ \frac{(t-3) \, (3 \, t-1) \, (3 \, t+1) \, (t^2+1)^2}{(t+2)^{10} \,
  (t+3)^{11}} =0. \] 
Hence $t$ can only be $3, \pm \frac{1}{3}$ or $\pm i$. Now
in fact $\Gamma$ is tri-involutive for each of these special values; for
instance, if $t = \frac{1}{3}$ then 
\[ \bA \ra 1, \quad \bB \ra \frac{1}{7}, \quad \bC \ra -\frac{1}{5}, \quad 
\bD \ra 0, \quad \bE \ra \frac{1}{3}, \quad \bF \ra -1 \] 
becomes an alignment. Replacing $L_2'$ by $M_3'$ leads
to the same inference by a similar route. In conclusion, the conditions $\CS_{123} =
\CS_{124} = \CS_{125} =0$ already force $\Gamma$ to be
tri-involutive. The main theorem is now completely proved. \qed 

\medskip 

In summary, we have an assembly of homomorphisms 
\[ G_\alpha \lra \SG(\F_\alpha), \quad \alpha \subseteq \six;  \] 
and the point of the theorem lies in the interconnections between
them. It is an intricate and highly regular structure  on the whole,
but at the moment I see no way of gaining any insight into it
except by a direct computational attack as above. It would be of interest to have a more
conceptual explanation for the entire phenomenon. 

\subsection{The undefined $g$-lines} \label{undefined.glines}
Some of the $g$-lines may become undefined on a
tri-involutive sextuple, i.e., the corresponding quadratic forms may
vanish identically. It is easy to classify all such cases by a
direct computation. Assume $\Gamma$ to be as in 
\S\ref{section.3inv.configuration}. 

\begin{itemize} 
\item 
The line $g(456)$ is undefined for any $p$; in fact all the Kirkman
points $K[w,456], w=1,2,3$ coincide with the pole of the line $Q_4Q_5Q_6$ in
Figure~\ref{fig:3inv}. 
\item
Additionally, if $p = \pm i$ then $g(\alpha)$ is undefined for every
$3$-element subset $\alpha \subseteq \{3,4,5,6\}$. 
\end{itemize} 

For a general $p$, the $g$-lines which remain defined are not all
distinct. We have 
\[ g(145)=g(245)=g(345) \] 
and similarly with $45$ replaced by $46$ and $56$. Thus there are only $13$ distinct lines. 

When $p = \pm i$, we have $g(134)=g(234)$ and similarly
with $34$ replaced by $35, 36, 45, 46, 56$. Thus there are only $10$
distinct lines. 

\subsection{} \label{section.chaslespoint}
We mention another characterisation of tri-involutivity. Given a triangle $\Delta = PQR$ in
$\P^2$, let $P'$ denote the pole of the line $QR$ with respect to
$\conic$, and similarly for $Q',R'$. Then $\Delta' = P'Q'R'$ is called
the polar triangle of $\Delta$. It is a
theorem due to Chasles (see~\cite[\S 7]{Coxeter}) that these two triangles are in perspective,
i.e., the lines $PP', QQ', RR'$ are concurrent, say in a point
$\tau(\Delta) = \tau(\Delta')$. It is also the case that the three
points 
\[ PQ \cap P'Q', \quad PR \cap P'R', \quad QR \cap Q'R' \] 
lie on the polar line of this point. In other words, the two triangles are in Desargues
configuration in such a way that the point of perspectivity and the axis of
perspectivity are in pole-polar relation. 

\begin{Proposition} \rm 
A sextuple $\Gamma$ is tri-involutive, if and only if it can be decomposed into two
triangles $\Delta_1, \Delta_2$ such that $\tau(\Delta_1) =
\tau(\Delta_2)$. 
\end{Proposition} 
\proof Assume $\Gamma$ to be tri-involutive. With $A, \dots, F$ as in~\S\ref{section.3inv.cr}, a
straightforward calculation shows that $\tau(ABC)$ and $\tau(DEF)$ are both
represented by the form 
\[ T = x_1^2 - x_1 \, x_2 + x_2^2. \] 
As to the converse, let $\Gamma = \Delta_1 \cup \Delta_2$ be a
decomposition as above. Since the 
construction of $\tau$ is compatible with automorphisms of $\conic$, 
we may assume that $\Delta_1$ has vertices $0, 1, \infty$ so that 
$\tau(\Delta_1)$ is given by $T$. If $\Delta_2$ has vertices $u, v, w$, then $\Delta_2'$ corresponds to 
\[ (x_1 - v \, x_2) \, (x_1 - w \, x_2), \quad 
(x_1 - u \, x_2) \, (x_1 - w \, x_2), \quad 
(x_1 - u \, x_2) \, (x_1 - v \, x_2). \] 
Since $\tau(\Delta_2')$ is also given by $T$, the three quadratic forms 
\[ (x_1 - u \, x_2)^2, \quad (x_1 - v \, x_2) \, (x_1 - w \, x_2),
\quad T \] 
must be linearly dependent and hence the $3 \times 3$ matrix of their
coefficients must have zero determinant. The same holds for the other two vertex
pairs. This reduces to a set of equations 
\[ \beta^2 - \alpha \, \beta - \alpha + 3 \, \beta +1 = u \, \alpha -
u - \alpha + \beta + 3 = u \, \beta + 1=0, \] 
where $\alpha = v+w$ and $\beta = v \, w$. The solutions are 
\[ \alpha = \frac{u^2-3 \, u+1}{u \, (u-1)}, \quad \beta = -
\frac{1}{u}, \] 
implying that
$\{ v, w \} = \{ \frac{u-1}{u}, \frac{1}{1- u} \}$. Hence $\Gamma = 
\{ 0, 1, \infty, u, \frac{u-1}{u}, \frac{1}{1-u} \}$ is tri-involutive. \qed 

\medskip 

The line $Q_4Q_5Q_6$ in Figure~\ref{fig:3inv} is represented by the form $T$
(see~\cite[\S 4.5]{CH1}), hence the point $\tau(ABC) = \tau(DEF)$ is
the pole of this line. It is not shown there explicitly, but in any
event must lie in the interior of the conic.

\section{The ideal of the tri-involutive locus}  

In this section we will determine the $SL_2$-equivariant minimal
resolution of the defining ideal of $\Omega$ using some
elimination-theoretic computations. The reader is referred to 
\cite{CH2} for a discussion of the invariant theory of binary
forms. (The set-up used there involves binary quintics rather than
sextics, but the general formalism is identical.) An exposition of the
same material may also be found in \cite[Ch.~4]{Sturmfels}. 

\subsection{} 
Let 
\[ \FF = \sum\limits_{i=0}^6 \, a_i \binom{6}{i} \, x_1^{6-i} \, x_2^i,
\]  denote the generic binary sextic with indeterminate coefficients,
and let 
\[ R = \complex[a_0,\dots,a_6] = \text{Sym}^\bullet \, (S_6), \] 
denote the coordinate ring of $\P^6$. Write 
\[ \Psi_p = x_1 \, x_2 \, (x_1 - \, x_2) \, (x_1 - p \, x_2) \, 
(x_1 - \frac{p-1}{p} \, x_2) \,  (x_1 - \frac{1}{1-p} \, x_2), \] 
for the sextic form representing a tri-involutive sextuple as in
\S\ref{section.3inv.configuration}. 

Let $J \subseteq R$ denote the homogeneous defining ideal of
$\Omega$; it consists of forms which vanish on the union of
$SL_2$-orbits of $\Psi_p$ taken over all $p$.  In order to calculate it, choose 
indeterminates $\alpha, \beta, \gamma, \delta$, and make substitutions 
\[ x_1 \ra \alpha \, y_1 + \beta \, y_2, \quad 
x_2 \ra \gamma \, y_1 + \delta \, y_2, \] 
into $(p-1) \, p \, \Psi_p$ to get a new sextic ${\mathfrak g}(y_1,y_2)$. Write it as 
\[ {\mathfrak g}(y_1, y_2) = \sum\limits_{i=0}^6 \;  \binom{6}{i} \, \varphi_i \,
y_1^{6-i} \, y_2^i, \] 
where $\varphi_i$ are polynomial expressions in the variables $\alpha, \beta,
\gamma, \delta, p$. This defines a ring homomorphism 
\[ \complex[a_0, \dots, a_6] \stackrel{g}{\lra} \complex[\alpha, \beta,
\gamma, \delta, p], \qquad a_i \lra \varphi_i, \] 
so that $J = \ker(g)$. I have calculated this ideal explicitly in 
{\sc Macaulay-2}; one pleasant surprise is that it turns out to be a
perfect ideal with minimal resolution
\[ 0 \la R/J \la R \la R(-5)^5 \la R(-6)^3 \oplus R(-7)
\la 0. \] 
In other words, $\Omega$ is an arithmetically Cohen-Macaulay
subvariety. Since $\Omega$ has Hilbert polynomial 
\[ \frac{2}{3} \, t^4 - \frac{5}{6} \, t^3 + \frac{35}{6} \, t^2 -
\frac{2}{3} \, t +2, \] 
its degree is $\frac{2}{3} \times 4! = 16$, in agreement with Proposition~\ref{prop.degree16}. 

We should like to identify the $SL_2$-representation corresponding to
the $5$-dimensional Betti module of minimal generators in degree $5$. It must be a
subrepresentation of $R_5 = \text{Sym}^5(S_6)$. By the
Cayley-Sylvester formula, 
\[ R_5 = S_2 \oplus S_2 \oplus S_4 \oplus \left\{ \text{summands of dimensions
  $>5$}\right\}, \] 
and hence on dimensional grounds the module can only be $S_4$. 
In other words, $J$ is minimally generated by the coefficients of a covariant of
degree-order $(5,4)$. The complete minimal system
for binary sextics is given in \cite[p.~156]{GY}; it shows that 
there is a unique such covariant up to a scalar. It can described as
follows: define 
\[ \vartheta_{24} = (\FF,\FF)_4, \quad 
\vartheta_{32} = (\vartheta_{24},\FF)_4, \quad 
\vartheta_{54} = (\vartheta_{24},\vartheta_{32})_1, \] 
where $\vartheta_{mq}$ is of degree-order $(m,q)$. Thus $J$ must be 
generated by the coefficients of $\vartheta_{54}$. This
implies that $\Gamma = \{z_1, \dots, z_6\}$ is
tri-involutive exactly when the covariant $\vartheta_{54}$ vanishes
identically on the sextic $\prod\limits_{i=1}^6 (x_1 - z_i \, x_2)$. 

The resolution shows that $J$ has a three-dimensional space of linear syzygies, and a unique quadratic
syzygy. In order to find the former, we must look for an identical
relation of the form $(\vartheta_{54}, \FF)_r=0$. Since binary sextics
has no covariant of degree-order $(6,2)$, perforce 
$(\vartheta_{54},\FF)_4 =0$. The quadratic syzygy is
similarly accounted for by the identity $(\vartheta_{54},\vartheta_{24})_4=0$,
since there is no invariant in degree $7$. In summary, the equivariant minimal resolution of $J$ is 
\[ 0 \la R/J \la R \la R(-5) \otimes S_4 \la R(-6) \otimes S_2 \oplus R(-7) \la 0. \] 

\subsection{} 
Since $\vartheta_{54}$ is defined to be the Jacobian of
$\vartheta_{24}$ and $\vartheta_{32}$, its vanishing implies a
functional dependency between the latter two. This is confirmed by a
direct computation. Indeed, 
\[ \vartheta_{24}(\Psi_p) = f_1(p) \, T^2, \quad \text{and} 
\quad \vartheta_{32}(\Psi_p) = f_2(p) \, T, \] 
where $T = x_1^2 - x_1 \, x_2 + x_2^2$ as in \S
\ref{section.chaslespoint}, and the $f_i$ are rational functions of
$p$. Thus, when $\FF$ is specialised to a tri-involutive
sextuple, $\vartheta_{32}$ evaluates to a quadratic form
which corresponds to the line containing its three centres of involution. 

\subsection{} \label{gorenstein3.case} 
A special case deserves to be mentioned. The precise
value of $f_1(p)$ is 
\[ f_1(p) = \frac{(p^2+1) \, (p^2-2 \, p + 2) \, (2 \, p^2 - 2 \, p +
  1)}{p^2 \, (p-1)^2}, \] 
which vanishes for 
\begin{equation} p = \pm i , \; 1 \pm i, \; \frac{1}{2} \, (1 \pm i). 
\label{special.p} \end{equation} 
Thus $\vartheta_{24}$ and $\vartheta_{32}$ are both identically zero for
these values of $p$, suggesting that the line containing the
three centres becomes `indeterminate' in some sense. As we will see, this is
indeed so. 

Specialise to $p=i$, i.e., let 
\[ A = 0, \quad B = 1, \quad C = \infty, \quad D = i, \quad E = \frac{i-1}{i} = 1
+i, \quad F = \frac{1}{1-i} = \frac{1}{2} (1+i). \] 
We have already seen in \S \ref{section.chaslespoint} that $\underbrace{\tau(ABC) = \tau(DEF)}_{T_1}$ for any
$p$. However, a direct calculation shows that in this case we additionally have 
\[ \underbrace{\tau(ADF) = \tau(BCE)}_{T_2}, \quad 
\underbrace{\tau(ACD) = \tau(BEF)}_{T_3}, \quad 
\underbrace{\tau(ABF) = \tau(CDE)}_{T_4}. \] 
Thus $\Gamma$ has altogether twelve centres of involution, which lie by
threes on four lines. 

The group $\RL(i)$ is strictly larger than the generic case. It is easily seen that $(\bA
\, \bE)$ and $(\bC \, \bF)$ are in $\RL(i)$; either element 
will turn all cross-ratios in~(\ref{crossratios.eq3}) into $-i$. Then
it follows automatically that $(\bB \, \bD) \in \RL(i)$. If $V$ denotes the group generated by
these two transpositions (so that $V$ is isomorphic to the Klein four-group),
then $\RL(i)$ is the internal direct product of $\RL$ and $V$. Hence $|\RL(i)| = 144$. 

The group $V$ permutes the four lines; for instance, 
transposing $A$ and $E$ will turn the pair of triangles $ABC, DEF$ into
$BCE, ADF$. The morphism 
$V \lra  \SG(\{T_1,T_2,T_3,T_4\})$ is given by 
\[ 
(\bA \, \bE) \lra (T_1 \, T_2) \, (T_3 \, T_4), \quad 
(\bC \, \bF) \lra (T_1 \, T_4) \, (T_2 \, T_3). 
\] 
\subsection{} 
In fact all the forms $\Psi_p$ for any of the six values
given in~(\ref{special.p}) lie in the same $SL_2$-orbit in $\P^6$, which we denote by
$\cZ$. This is a well-known geometric object. The orbit is Zariski closed in $\P^6$;
indeed there are very few orbits of binary forms which have this property, and they have all
been classified in \cite{AF}. (There it is described as the orbit of $x_1^5 \, x_2  - x_1
\, x_2^5$, which comes to the same thing.) Moreover, $\cZ \subseteq \P^6$ is an arithmetically
Gorenstein subvariety in codimension $3$, and its ideal $I$ is generated by the coefficients of
$\vartheta_{24}$ (see~\cite[\S 3]{MU}). It has a self-dual Buchsbaum-Eisenbud resolution 
\[ 0 \la R/I \la R \la R(-2) \otimes S_4 \la R(-3) \otimes S_4 \la
R(-5) \la 0. \] 

Unfortunately one cannot draw a diagram of this sextuple and its twelve centres of involution, since the 
equality $\langle A, B, C, F \rangle = i$ implies that not all six points can be chosen to be simultaneously
real.

\medskip 

\centerline{--} 

\vspace{1cm}

\parbox{7cm}{ \small 
Jaydeep Chipalkatti \\
Department of Mathematics \\ 
University of Manitoba \\ 
Winnipeg, MB R3T 2N2 \\ 
Canada. \\ \\
{\tt chipalka@cc.umanitoba.ca}}

\end{document}